\documentclass{amsart}
\usepackage{amsmath}
\usepackage{amssymb}
\usepackage{kotex}
\usepackage{enumerate}
\usepackage{amsbsy}
\usepackage{amsfonts}
\usepackage{color}
\usepackage{parskip}

\newtheorem{thm}{Theorem}[section]
\newtheorem{lem}[thm]{Lemma}
\newtheorem{prop}[thm]{Proposition}

\newtheorem{rem}[thm]{Remark}

\numberwithin{equation}{section}

\begin{document}

	\title[4D Maxwell-Klein-Gordon System in the Lorenz gauge]{ Global well-posedness and scattering of the energy-critical Maxwell-Klein-Gordon system in the Lorenz gauge}

    \author[S. Hong]{Seokchang Hong}
    \address{Department of Mathematical Sciences, Seoul National University, Seoul 08826, Republic of Korea}
    \email{seokchangh11@snu.ac.kr}

	
	\begin{abstract}
	We study initial value problem of the $(1+4)$-dimensional Maxwell-Klein-Gordon system (MKG) in the Lorenz gauge. 
	Since (MKG) in the Lorenz gauge does not possess an obvious null structure, it is not easy to handle the nonlinearity. 
	To overcome this obstacle, we impose an additional angular regularity. 
	In this paper, we prove global well-posedness and scattering of (MKG) for small data in a scale-invariant space which has extra weighted regularity in the angular variables. Our main improvement is to attain the scaling critical regularity exponent and prove global existence of solutions to (MKG) in the Lorenz gauge. 
	\end{abstract}

		\maketitle

\section{Introduction}

In this paper, we investigate global well-posedness and scattering of the $(1+4)$-dimensional Maxwell-Klein-Gordon system in the Lorenz gauge. The Maxwell-Klein-Gordon (MKG) system describes a physical phenomena of a spin-zero particle in an electromagnetic field.
The (MKG) system is obtained by coupling the Klein-Gordon scalar field $\phi:\mathbb R^{1+4}\rightarrow\mathbb C$ with an elctromagnetic field $F$. To be precise, we consider the covariant form of the (MKG) system:
\begin{align}\label{mkd}
\begin{aligned}
\partial^\nu F_{\mu\nu}  & = \textrm{Im}(\phi\overline{\mathcal D_\mu\phi}), \\
\mathcal D^\mu\mathcal D_\mu A & = m^2\phi,
\end{aligned}
\end{align}
where $F=dA$ is the associated curvature $2$-form given by $F_{\mu\nu}=\partial_\mu A_\nu-\partial_\nu A_\mu$ and $A=(A_\mu),\ \mu=0,1,\cdots,4$ is a real-valued $1$-form on $\mathbb R^{1+4}$. We write $\mathcal D^\mu=\partial^\mu-iA^\mu$ for the connection. Here the constant $m\ge0$ is a mass. Since we study the initial value problem of (MKG) in a scale-invariant Besov space, we put $m=0$ in the sequel. In other words, we consider the mass-less (MKG) system. 

The (MKG) in $\mathbb R^{1+4}$ obeys the law of conservation of energy. The conserved energy of a solution $(A,\phi)$ at time $t$ is defined as  
\begin{align}
\mathcal E_{MKG}(\phi,A):=\frac12\int_{\mathbb R^4}|F_{\mu\nu}(t)|^2+|\mathcal D_\mu \phi(t)|^2\,dx. 
\end{align}
We also note that the (MKG) system is invariant under the scaling 
$
(A,\phi)\mapsto (\lambda^{-1}A,\lambda^{-1}\phi)(\lambda^{-1}t,\lambda^{-1}x),
$
and hence the critical Sobolev space for the system \eqref{mkd} is $\dot{H}^1$. Consequently, the Cauchy problem of the $(1+4)$-dimensional (MKG) system can be regarded as energy-critical problem.

We remark that the gauge potential $A_\mu$ is not necessarily unique. Indeed, for a sufficiently smooth real-valued function $\Lambda$ on $\mathbb R^{1+4}$, the (MKG) system is invariant under the transform $(A,\phi)\mapsto(A-d\Lambda,e^{i\Lambda}\phi)$. 
Hence we have gauge freedom and it is reasonable to choose a representative, which is suitable for our purpose. (See also \cite{selbtes}.)  

\subsection{Maxwell-Klein-Gordon equations in the Lorenz gauge}
Imposing the Lorenz gauge $\partial^\mu A_\mu=0$ we obtain the following nonlinear wave equations after some computation
\begin{align}\label{mkd-lorenz}
\begin{aligned}
\Box \phi = 2iA^\mu\partial_\mu\phi+A^\mu A_\mu\phi, \\
\Box A = -\textrm{Im}(\phi\overline{\partial\phi})-A|\phi|^2,
\end{aligned}
\end{align} 
where $\Box=\partial_t^2-\Delta$ is the d'Alembertian.
The Maxwell-Klein-Gordon system has been extensively studied \cite{klaima1,selbtes}. Typical choices of gauge are the Lorenz gauge $\partial^\mu A_\mu=0$ and the Coulomb gauge $\partial^jA_j=0$. 
 Keel-Roy-Tao  \cite{keelroytao} proved global existence below the energy space and almost optimal well-posedness is proven by Machedon-Sterbenz \cite{masterbenz}. In the Lorenz gauge, Pecher \cite{pecher1} showed almost critcal well-posedness in the Fourier-Lebesgue spaces. 
 
 From now on we exclusively consider $(1+4)$-dimensional setting. The (MKG) system is relatively well-understood in the Coulomb gauge. Selberg \cite{selberg} proved well-posedness in $H^{1+}$, which is almost optimal and Krieger-Sterbenz-Tataru \cite{kriegersterbenztataru} showed global well-posedness with small $H^1$ norm. Then Oh-Tataru \cite{ohtataru,ohtataru1,ohtataru2} obtained global well-posedness with arbitrarily large data. 
 Concerning the Lorenz gauge, Pecher \cite{pecher} proved local well-posedness for $H^{\frac76+}$ data. Thus well-posedness for critical regularity data in the Lorenz gauge is still open. 
 
 We briefly discuss the technical difficulty in proving low regularity well-posedness. When we are concerned with quadratic nonlinearity, the most difficult type of interaction is when the two inputs give rise to an output which is close to the light cone in the space-time Fourier space.   However, if nonlinearity has a certain cancellation property, it would be possible to obtain better regularity properties. 
  Such cancellation typically given by null structure plays a crucial role to attain the critical regularity. We refer the readers to  \cite{klaima,fosklai,klaitataru,leevargas} for the study on the null form estimates and their applications. 
   
   In the Lorenz gauge, however, the main obstacle is that the null structure inside the equations \eqref{mkd-lorenz} is \textit{not enough} to attain the critical regularity as we are concerned with a low dimensional setting. 
   By \textit{not enough}, we mean that the cancellation property of parallel interaction given by the null structure is only the angle $\angle(\xi_1,\xi_2)$ between the two input frequencies.  
   If the null structure is \textit{enough}, for instance, it yields more than $\angle(\xi_1,\xi_2)$, such as the type of null form $Q_0$, defined by $Q_0(u,v)=\partial_tu\partial_tv-\partial_ju\partial^jv$ which gives $\angle(\xi_1,\xi_2)^2$, we could overcome this problem. Unfortunately, we cannot expect such an enough null structure in this case
    and hence it is not easy to handle the parallel interactions. 
    
   Inspired by the work of Sterbenz \cite{sterbenz1}, we expect that the rotation generators $\Omega_{ij}=x_i\partial_j-x_j\partial_i$ plays a crucial role since they eliminate such delicate interactions. 
 Furthermore, we will get a significant gain over the \textit{classical} Strichartz estimates for the wave equation \cite{sterbenz1,cholee}, and hence obtain a crucial improvement at the level of multilinear estimates.  
 See also \cite{candyherr,candyherr1,sterbenz2,wang} and references therein for the study on the Dirac-Klein-Gordon and Yang-Mills systems via angular regularity. 
 In this paper, we study global well-posedness and scattering of the mass-less (MKG) in a scale-invariant Besov space with additional angular regularity. For the purpose, we consider the following function space
\begin{align}
\|f\|_{\dot{B}^{s,\sigma}_\Omega} := \|\langle\Omega\rangle^\sigma f\|_{\dot{B}^s}, 
\end{align}
where $\langle\Omega\rangle^\sigma  :=   (1-\Delta_{\mathbb S^3})^\frac\sigma2$.
Here, $\dot{B}^s=\dot{B}^s_{2,1}$ is the usual homogeneous Besov space and $\Delta_{\mathbb S^3}$ is the Laplace-Beltrami operator on the unit sphere $\mathbb S^3$. 
We consider our system \eqref{mkd-lorenz} with the initial data
\begin{align*}
(\phi,\partial_t\phi)|_{t=0}&:=(\phi_0,\phi_1)\in \dot{B}^{1,1}_\Omega\times\dot{B}^{0,1}_\Omega, \\ 
(A_0,\mathbf A)|_{t=0}&:=(a_0,\mathbf a)\in\dot{B}^{1,1}_\Omega.
\end{align*} 
We also write $\partial_tA_\mu|_{t=0}=\dot{a_\mu}\in\dot{B}^{0,1}_\Omega$. Then we consider the following constraints for data:
\begin{align}
a_0 = \dot{a_0} = 0, \label{ini1} \\
\partial^ka_k =0	. \label{ini3}
\end{align}
We need the condition \eqref{ini1} because otherwise the Lorenz gauge condition does not determine the potential uniquely. The condition \eqref{ini3} follows from the Lorenz condition in connection with \eqref{ini1}. (See also \cite{pecher}.) We state our main result.



\begin{thm}[Global well-posedness]\label{gwp}
Suppose that the initial data $(\phi_0,\phi_1,\mathbf a,\dot{\mathbf a})$ satisfy the smallness condition
\begin{align}\label{smallness}
\|(\phi_0,\phi_1,\mathbf a,\dot{\mathbf a})\|_{\dot{B}^{1,1}_\Omega\times\dot{B}^{0,1}_\Omega\times\dot{B}^{1,1}_\Omega\times\dot{B}^{0,1}_\Omega  } \le \epsilon_0,
\end{align}
and the constraints \eqref{ini1} and \eqref{ini3}. Then
there exists a global solution $(\phi(t),A(t))$ to the mass-less (MKG) system with the Lorenz gauge, which satisfies the following stability condition:
\begin{align}\label{stability}
\sup_{-\infty<t<+\infty}\|(\phi(t),\partial_t\phi(t), A(t),\partial_tA(t))\|_{\dot{B}^{1,1}_\Omega\times\dot{B}^{0,1}_\Omega\times\dot{B}^{1,1}_\Omega\times\dot{B}^{0,1}_\Omega } \le C.
\end{align}
Furthermore, the solution depends smoothly on the initial data.
\end{thm}
Our main improvement is to attain the scaling critical regularity exponent and prove global existence of solutions to (MKG) using a scale-invariance of function space. Then scattering is followed by Theorem \ref{gwp}.
\begin{thm}[Scattering]\label{scattering}
	For any given initial data $(\phi_0,\phi_1,\mathbf a,\dot{\mathbf a})\in \dot{B}^{1,1}_\Omega\times\dot{B}^{0,1}_\Omega\times\dot{B}^{1,1}_\Omega\times\dot{B}^{0,1}_\Omega $ satisfying the conditions \eqref{ini1}, \eqref{ini3}, and \eqref{smallness}, there exist unique functions $(\phi_0^+,\phi_1^+,\mathbf a^+,\dot{\mathbf a}^+)$ and $(\phi_0^-,\phi_1^-,\mathbf a^-,\dot{\mathbf a}^-)$ such that 
	\begin{align}
	\begin{aligned}
	\lim_{t\rightarrow\pm\infty}\big( &\|\phi(t)-\phi^\pm(t)\|_{\dot{B}^{1,1}_\Omega}+\|\partial_t\phi(t)-\partial_t\phi^\pm(t)\|_{\dot{B}^{0,1}_\Omega} \\
	&\qquad\qquad\qquad\qquad+\|A(t)-A^\pm(t)\|_{\dot{B}^{1,1}_\Omega}+\|\partial_tA(t)-\partial_tA^\pm(t)\|_{\dot{B}^{0,1}_\Omega}  \big) = 0.	
	\end{aligned}
	\end{align}
Moreover, the scattering operator which maps $(\phi_0,\phi_1,\mathbf a,\dot{\mathbf a})$ to $(\phi_0^\pm,\phi_1^\pm,\mathbf a^\pm,\dot{\mathbf a}^\pm)$ is a local diffeomorphism in $\dot{B}^{1,1}_\Omega\times\dot{B}^{0,1}_\Omega\times\dot{B}^{1,1}_\Omega\times\dot{B}^{0,1}_\Omega $.
\end{thm}

\subsection{Strategy of proof}
We follow the approach due to the works of Sterbenz and Wang \cite{sterbenz2,wang}. That is, we shall define the appropriate function spaces and use bilinear decomposition for angles and estimate the nonlinearity in the function spaces. As the author of \cite{sterbenz2} mentioned, the estimate of cubic terms is very straightforward. Indeed, we will only use H\"{o}lder's inequality and Strichartz estimates, and hence estimates of bilinear forms in \eqref{mkd-lorenz} will be the crucial part of this paper.
To do this, we consider all possible frequency interactions such as High$\times$High and Low$\times$High interactions. 
By the mercy of an extra weighted regularity in angular variables, we enjoy the improved space-time Strichartz estimates and hence the High$\times$High interaction can be treated rather easier than the Low$\times$High case. (See Section \ref{pf-hh-near-cone}.)

We remark that in the Low$\times$High interaction, the situation that the low frequency controls the angular regularity becomes more difficult case, since we cannot exploit the angular concentration estimate. 
(See Remark \ref{rem-lh-ang}.)

\noindent\textbf{Organisation.} The rest of this paper is organised as follows. In Section 2, we recall the Strichartz estimates and reveal null structure of $A^\mu\partial_\mu\phi$. We construct the function spaces via the Littlewood-Paley projections in Section 3. We introduce the angular decompositions of bilinear form in Section 4. Then Section 5 will be the main part of this paper, devoted to the proof of our main result.

\noindent\textbf{Notations.} 
As usual different positive constants, which are independent of dyadic numbers $\mu,\lambda$, and $d$ are denoted by the same letter $C$, if not specified. $A \lesssim B$ and $A \gtrsim B$ means that $A \le CB$ and
$A \ge C^{-1}B$, respectively for some $C>0$. $A \approx B$ means that $A \lesssim B$ and $A \gtrsim B$, i.e., $\frac1CB \le A\le CB $ for some absolute constant $C$. We also use the notation $A\ll B$ if $A\le \frac1CB$ for some large constant $C$. Thus for quantites $A$ and $B$, we can consider three cases: $A\approx B$, $A\ll B$ and $A\gg B$. In fact, $A\lesssim B$ means that $A\approx B$ or $A\ll B$. We shall use the notation $A\pm$ which means that for small positive $\epsilon>0$, we may replace $A\pm$ by $A\pm\epsilon$. For example, we shall write the improved Strichartz estimates with additional angular regularity as
$$
\|e^{\mp it|\nabla|}f_1\|_{L^2_tL^{3+}_x} \lesssim \|\langle\Omega\rangle^\frac12 f_1\|_{L^2_x}.
$$
The spatial and space-time Fourier transforms are defined by 
\begin{align*}
\widehat{f}(\xi) &= \int_{\mathbb R^4} e^{-ix\cdot\xi}f(x)\,dx, \\ 
\widetilde{u}(\tau,\xi) &= \int_{\mathbb R^{1+4}}e^{-i(t\tau+x\cdot\xi)}u(t,x)\,dtdx.
\end{align*}
We also write $\mathcal F_x(f)=\widehat{f}$ and $\mathcal F_{t,x}(u)=\widetilde{u}$. We denote the backward and forward wave propagation of a functiom $f$ on $\mathbb R^4$ by
$$
e^{\mp it|\nabla|}f = \int_{\mathbb R^4}e^{ix\cdot\xi}e^{\mp it|\xi|}\widehat{f}(\xi)\,d\xi.
$$
The notation $V$ denotes  the parametrix for the inhomogeneous wave equation with zero initial data, that is, $u=Vf$ if and only if
$$
\Box u=f, \quad u(0,\cdot) = 0, \quad \partial_tu(0,\cdot) = 0. 
$$
Let $E$ be any fundamental solution to the homogeneous wave equation, that is, $\Box E=\delta$, where $\delta$ is the Dirac delta distribution. Then we can represent the parametrix operator $V$ via the following formula:
$$
V(f) = E*f-W(E*f),
$$
where, for any smooth well-defined function $g=g(t,x)$, $W(g)$ denotes the solution of linear homogeneous wave equation with initial data $(g,\partial_tg)$.

Finally, we recall some basic facts from harmonic analysis on the sphere. The most of the ingredients here can be found in \cite{sterbenz2}. We also refer the readers to \cite{steinweiss} for the systematic introduction to the analysis on the sphere. We denote the standard infinitesimal generators of the rotations on $\mathbb R^4$ by
$\Omega_{ij}=x_i\partial_j-x_j\partial_i$. Then the Laplace-Beltrami operator $\Delta_{\mathbb S^3}$ can be written as
$$
\Delta_{\mathbb S^3} = \sum_{i<j}\Omega_{ij}^2.
$$ 
We define the fractional order of the spherical Laplacian
$$
 |\Omega|^\sigma = (-\Delta_{\mathbb S^3})^\frac\sigma2.
$$ 
We also use the inhomogeneous form of $\Delta_{\mathbb S^3}$:
$$
\langle\Omega\rangle^\sigma f = f_0+|\Omega|^\sigma f,
$$
where $f_0$ is the radial part of $f$, given by
$$
f_0(r) = \frac{1}{|\mathbb S^3|}\int_{\mathbb S^3}f(r\omega)\,d\omega.
$$
An important fact of the operators $\langle\Omega\rangle^\sigma$ is that they are commutative with the Fourier transform:
$$
\mathcal F[\langle\Omega\rangle^\sigma f] = \langle\Omega\rangle^\sigma \mathcal F(f). 
$$
Then the homogeneous Besov space with additional angular regularity is defined as
$$
\|f\|_{\dot{B}^{s,\sigma}_\Omega} := \|\langle\Omega\rangle^\sigma f\|_{\dot{B}^s},
$$
where $\|f\|_{\dot{B}^s}=\sum_{\lambda\in 2^{\mathbb Z}}\lambda^s\|P_\lambda f\|_{L^2_x}$, and $P_\lambda$ is the Littlewood-Paley prejection onto the set $\{ \xi\in\mathbb R^4 : |\xi|\approx\lambda\}$.
\begin{rem}
Strictly speaking, the following Leibniz rule is not true:
$$
\langle\Omega\rangle(f g) = (\langle\Omega\rangle f)g+f(\langle\Omega\rangle g).
$$
This is because obviously the operator $|\Omega|$ is a non-local operator. On the other hand, $\Omega_{ij}$ is clearly a local operator and the Leibniz rule holds. However, for convenience, by abuse of notation we replace any instance of an single $\Omega_{ij}$ with the operator $\langle\Omega\rangle$ and assume the above Leibniz rule for $\langle\Omega\rangle$ is true. 
See also \cite{sterbenz2, wang}. 
\end{rem}


\section{Strichartz estimates and Null structure}
\subsection{Strichartz estimates}
We first introduce the classical Strichartz estimates for the homogeneous wave equations, based soley on translation invariant derivatives of the initial data. (See also \cite{strichartz,keeltao}.) 
\begin{prop}
Let $n=4$ be the number of spatial dimensions, and let $\sigma=\frac32$ be the corresponding Strichartz admissibility exponent. If $f$ is any function of the spatial variable only, denote by $f_1=P_1f$ its unit frequency projection. Then one has the following family estimates for $2\le q$:
\begin{align}
\|e^{\mp it|\nabla|}f_1\|_{L^q_tL^r_x} \lesssim \|f_1\|_{L^2_x},
\end{align}
where $\frac1q+\frac\sigma r\le\frac\sigma2$.
\end{prop}
In sequel, we only use a few selected subset of admissible pair $(q,r)$, namely,
\begin{align}
\|e^{\mp it|\nabla|}f_1\|_{L^\infty_tL^2_x} & \lesssim \|f_1\|_{L^2_x}, \\
\|e^{\mp it|\nabla|}f_1\|_{L^2_tL^\infty_x} & \lesssim \|f_1\|_{L^2_x}, \\
\|e^{\mp it|\nabla|}f_1\|_{L^2_tL^6_x} & \lesssim \|f_1\|_{L^2_x}.
\end{align}
The sharpness of the above space-time estimates is proven by Knapp type counterexample, which does not have radial symmetry. Hence it is natural to expect that the Strichartz estimates would be improved if one impose spherical symmetry. 
We refer to \cite{sterbenz1,cholee} for the details and other differential operators.  
\begin{prop}
Let $n=4$ be the number of spatial dimensions, and let $\sigma_\Omega=3$ denote the four dimensional angular Strichartz admissible exponent. Let $f_1$ be a unit frequency function of the spatial variable only. Then for indices $(q,r)$ such that $\frac1q+\frac\sigma r\ge\frac\sigma2$ and $\frac1q+\frac{\sigma_\Omega}{r}<\frac{\sigma_\Omega}{2}$, and for every $0<\epsilon$, there is a $C_\epsilon$ which depends only on $\epsilon$ such that the following estimates hold:
\begin{align}
\|e^{\mp it|\nabla|}f_1\|_{L^q_tL^r_x} \lesssim C_\epsilon\|\langle\Omega\rangle^sf_1\|_{L^2_x},	
\end{align}
	where $s=(1+\epsilon)(\frac3r+\frac2q-\frac32)$.
\end{prop}
In particular, we shall use often the following space-time estimate:
\begin{align}
	\|e^{\mp it|\nabla|}f_1\|_{L^2_tL^{3+}_x} &\lesssim \|\langle\Omega\rangle^\frac12 f_1\|_{L^2_x}.
\end{align}

\subsection{Null form of $A^\mu\partial_\mu\phi$}
In this section, we derive the null form in $A^\mu\partial_\mu\phi$. See also \cite{pecher}. We begin with the standard $Q$-type null form introduced in \cite{klaima}:
\begin{align*}
Q_{0j}(u,v) & = \partial_tu\partial_jv-\partial_ju\partial_tv, \\
Q_{jk}(u,v)& = \partial_ju\partial_kv-\partial_ku\partial_jv.
\end{align*}
We decompose the spatial part $\mathbf A$ of the gauge potential $A$ into the divergence-free and curl-free parts
\begin{align*}
\mathbf A = \mathbf A^{\rm df}+\mathbf A^{\rm cf},
\end{align*}
where
\begin{align*}
A^{\rm df}_j := R^k(R_jA_k-R_kA_j),\quad A^{\rm cf}_j := -R_jR_kA^k.
\end{align*}
Here $R_j=|\nabla|^{-1}\partial_j$ is the Riesz transform. Then we see that
\begin{align*}
A^\mu\partial_\mu\phi &= A^0\partial_t\phi+\mathbf A\cdot\nabla\phi	 \\
& = -A_0\partial_t\phi+\mathbf A^{\rm cf}\cdot\nabla\phi+\mathbf A^{\rm df}\cdot\nabla \phi \\
& =: \mathcal N_1+\mathcal N_2.
\end{align*}
By the Lorenz gauge condition $\partial_tA_0=\partial^jA_j$, we get
\begin{align*}
\mathcal N_1 & = -A_0\partial_t\phi - R_jR_kA^k\partial^j\phi \\
& = -A_0\partial_t\phi - (|\nabla|^{-1}\partial_tA_0)\cdot\nabla\phi \\
& = \partial_j(|\nabla|^{-1}R^jA_0)\partial_t\phi - \partial_t(|\nabla|^{-1}R_jA_0)\partial^j\phi	\\
& = -Q_{0j}(|\nabla|^{-1}R^jA_0,\phi).
\end{align*}
We also have
\begin{align*}
\mathcal N_2 & = R^k(R_jA_k-R_kA_j)\partial^j\phi \\
& = \langle\nabla\rangle^{-2}	\partial^k\partial_jA_k\partial^j\phi +A_j\partial^j\phi \\
& = -\frac12\left( \langle\nabla\rangle^{-2}(\partial_j\partial^jA_k-\partial_j\partial_kA^j)\partial^k\phi - \langle\nabla\rangle^{-2}(\partial^k\partial_jA_k-\partial^k\partial_kA_j)\partial^j\phi \right) \\
& = -\frac12\left( \partial_j\langle\nabla\rangle^{-1}(R^jA_k-R_kA^j)\partial^k\phi - \partial^k\langle\nabla\rangle^{-1}(R_jA_k-R_kA_j)\partial^j\phi \right) \\
& = -\frac12 Q_{jk}(\langle\nabla\rangle^{-1}(R^jA^k-R^kA^j),\phi).
\end{align*}
We shall denote the Fourier symbols of $Q_{0j}$ and $Q_{jk}$ by $q_{0j}(\xi,\eta)$ and $q_{jk}(\xi,\eta)$, respectively. Then the symbols satisfy the following estimates:
\begin{align}
|q_{0j}(\xi,\eta)|,|q_{jk}(\xi,\eta)| \lesssim |\xi||\eta|\angle(\xi,\eta),
\end{align}
where $\angle(\xi,\eta)=\arccos(\frac{\xi\cdot\eta}{|\xi||\eta|})$ is the angle between $\xi$ and $\eta$. (See \cite{selbtes1}.) \\
The above notations seem too complicated. 
However, an important point here is that the bilinear form $A^\mu\partial_\mu\phi$ contributes an additional angle between two input frequencies. In other words, it is not necessary to distinguish the null forms $Q_{0j}$ and $Q_{jk}$. In the sequel, we will by abuse of notation simply write $Q(\varphi,\phi)$ for the bilinear form $A^\mu\partial_\mu\phi$. That is, we even ignore the vector components of $A^\mu$ and denote it by $\varphi$ shortly. 

\section{Function spaces}

This section is devoted to the introduction of preliminary setup to be used in the proof of Theorem \ref{gwp}. We shall establish the function spaces in this section. The readers can find the following notations in \cite{sterbenz2,wang}. 
\subsection{Multipliers}
Let $\beta$ be a smooth bump function given by $\beta(s)=1$ for $|s|\le1$ and $\beta(s)=0$ for $|s|\ge2$. We define the dyadic scaling of $\beta$ by $\beta_\lambda(s)=\beta(\lambda^{-1}s)$ for $\lambda\in2^{\mathbb Z}$. Now we define the Fourier multipliers which give localisation with respect to the spatial frequency, space-time frequency, distance to the cone (modulation), distance to the lower cone and the upper cone as follows:
\begin{align}
p_\lambda(\xi) & = \beta_{2\lambda}(|\xi|)-\beta_{\lambda/2}(|\xi|), \quad s_\lambda(\tau,\xi) = \beta_{2\lambda}(|(\tau,\xi)|)-\beta_{\lambda/2}(|(\tau,\xi)|), \\
c_d(\tau,\xi) & = \beta_{2d}(|\tau|-|\xi|)-\beta_{d/2}(|\tau|-|\xi|), \\
c_d^+(\tau,\xi) & = \beta_{2d}(\tau+|\xi|)-\beta_{d/2}(\tau+|\xi|),\quad c_d^-(\tau,\xi) = \beta_{2d}(\tau-|\xi|)-\beta_{d/2}(\tau-|\xi|).
\end{align}
Then we define the corresponding Fourier projection operators. For example, $\mathcal F_x(P_\lambda f)=p_\lambda\mathcal F_x(f)$ and $\mathcal F_{t,x}(S_\lambda u)=s_\lambda\mathcal F_{t,x}(u)$. We also write $S_{\lambda,d}=S_\lambda C_d$, $S^\pm_{\lambda,d}=S_\lambda C^\pm_d$, and denote
$$
S_{\lambda,\cdot\le d} = \sum_{\delta\le d}S_{\lambda,\delta},\quad S^\pm_{\lambda,\cdot\le d} = \sum_{\delta\le d}S^\pm_{\lambda,\delta},
$$
and the projection $S_{\lambda,d\le\cdot}$ is defined in the obvious way. 

For any given $0<\eta\lesssim1$, we define $\mathcal C_\eta$ to be a collection of finitely overlapping caps of size $\eta$ in the unit sphere $\mathbb S^3$ in $\mathbb R^4$. Then we denote a smooth partition of unity subordinate to the angular sectors $\{ \xi\neq0 : \frac{\xi}{|\xi|}\in\omega \}$ by $\{b^\omega_\eta\}_{\omega\in\mathcal C_\eta}$. The corresponding angular localisation operator is denoted by $B^\omega_\eta$. For simplicity, we allow the abuse of notation and identify the angular sectors $\omega\in\mathcal C_\eta$ and their centre. 
We define
\begin{equation}
S^\omega_{\lambda,d} = B^\omega_{(\frac d\lambda)^\frac12}P_\lambda S_{\lambda,d}, \quad S^\omega_{\lambda,\cdot\le d} =  B^\omega_{(\frac d\lambda)^\frac12}P_\lambda S_{\lambda,\cdot\le d}. 
\end{equation}
We have the following lemma concerning boundedness of the above localisation operators. (See also \cite{sterbenz2,wang}.) 
\begin{lem}\label{multi-bdd}
\begin{enumerate}
\item The following multipliers are given by $L^1_{t}L^1_x$ kernels and are uniformly bounded in $L^1_{t}L^1_x$: $\lambda^{-1}\nabla S_\lambda,\ B^\omega_{(\frac d\lambda)^\frac12}P_\lambda, \ S^\omega_{\lambda,d}, \ (\lambda d)VS^\omega_{\lambda,d}$ and also those operators are bounded in mixed Lebesgue spaces $L^q_tL^r_x$. \\
\item The following multipliers are uniformly bounded in the $L^q_tL^2_x$ spaces for $1\le q\le+\infty$: $S_{\lambda,d}, \ S_{\lambda,\cdot\le d}$. 
\end{enumerate}
\end{lem}
\subsection{Function spaces}
Now we establish the function spaces with space-time frequency localised in the support of $s_\lambda(\tau,\xi)$. We let
$$
\|u\|_{X^{\frac12,p}_\lambda}^p := \sum_{d;d\in2^\mathbb Z} d^\frac p2 \|S_{\lambda,d}u\|_{L^2_{t}L^2_x}^p,
$$
and
$$
\|u\|_{Y_\lambda} := \lambda^{-1}\|\Box S_\lambda u\|_{L^1_tL^2_x}.
$$
We also define the space $Z_{\Omega,\lambda}$ by
$$
\|f\|_{Z_{\Omega,\lambda}} := \lambda^{-1}\sum_{d\lesssim\lambda}\int\sup_\omega\|S^\omega_{\lambda,d}f\|_{L^\infty_x}\,dt.
$$
Then we have
\begin{align*}
\|f\|_{Z_{\Omega,\lambda}} \lesssim  \lambda^{-1}\|S_\lambda\Box\langle\Omega\rangle f\|_{L^1_tL^2_x} \lesssim \|S_\lambda\langle\Omega\rangle f\|_{Y_\lambda},
\end{align*}
and hence $\langle\Omega\rangle^{-1}Y_\lambda\subset Z_{\Omega,\lambda}$. We also have
\begin{align}\label{prop-z}
\sup_{d\lesssim\lambda}\lambda^{-1}\left\|\sup_\omega \|S^\omega_{\lambda,d}f\|_{L^\infty_x}\right\|_{L^1_t} \lesssim \|f\|_{F_{\Omega,\lambda}}.
\end{align}
 Finally, we define the function space
$$
F_{\Omega,\lambda} = \langle\Omega\rangle^{-1}(X^{\frac12,1}_\lambda+Y_\lambda)\cap S_\lambda(L^\infty_tL^2_x)\cap Z_{\Omega,\lambda}.
$$
As $X^{\frac12,1}_\lambda\subset S_\lambda(L^\infty_tL^2_x)$, naturally we have $\Box X^{\frac12,1}_\lambda\subset\Box S_\lambda(L^\infty_tL^2_x)$, and by duality, we obtain
\begin{align}
VS_\lambda(L^1_tL^2_x)\subset (\Box X^{\frac12,1}_\lambda)' = (\lambda X^{-\frac12,1}_\lambda)' = \lambda^{-1}X^{\frac12,\infty}_\lambda,
\end{align}
and recall an obvious relation $X^{\frac12,1}_\lambda\subset X^{\frac12,\infty}_\lambda$, hence we have
\begin{align}
d^\frac12 \|S_{\lambda,d}f\|_{L^2_{t,x}} & \lesssim \|\langle\Omega\rangle^{-1}f\|_{F_{\Omega,\lambda}}, \quad \textrm{ for } d\in 2^{\mathbb Z}, \quad 0<d\le\lambda, \\
d^\frac12 \|S_{\lambda,d\le\cdot}f\|_{L^2_{t,x}} & \lesssim \|\langle\Omega\rangle^{-1}f\|_{F_{\Omega,\lambda}}.
\end{align}
Now we define the Besov type function space, which we will iterate to prove our main theorem:
\begin{align}
\|u\|_{F_\Omega} := \sum_{\lambda\in2^{\mathbb Z}}\lambda\|S_\lambda u\|_{F_{\Omega,\lambda}}.	
\end{align}
We list several space-time estimates to be used throughout this paper.
\begin{align}
\|S_\mu\langle\Omega\rangle u\|_{L^\infty_tL^2_x} & \lesssim \|u\|_{F_{\Omega,\mu}}, \\
\|S_\mu\langle\Omega\rangle u\|_{L^2_tL^\infty_x} & \lesssim \|u\|_{F_{\Omega,\mu}}, \\
\|S_\mu\langle\Omega\rangle u\|_{L^2_tL^6_x} & \lesssim \mu^\frac56\|u\|_{F_{\Omega,\mu}}, \\
\|S_\mu\langle\Omega\rangle^\frac12u\|_{L^2_tL^{3+}_x} & \lesssim \mu^{\frac16+}\|u\|_{F_{\Omega,\mu}}	.
\end{align}
We end this section with angular concentration estimates and Sobolev embedding estimates. 
\begin{lem}[Lemma 5.2 of \cite{sterbenz2}]\label{ang-con-est}
Let $2\le n$ be a given integer. Then for any test function $u$ on $\mathbb R^n$, and any $2\le p<\infty$, one has the following estimate:
$$
\sup_\omega\|B^\omega_\eta u\|_{L^p} \lesssim \eta^s\|\langle\Omega\rangle^su\|_{L^p},
$$
where $0\le s<\dfrac{n-1}{p}$.
\end{lem} 
\begin{lem}\label{sobolev-emb}
Let $f$ be a test function on $\mathbb R^4$. Then one has the following frequency localised estimate:
\begin{align}
\|B^\omega_\eta P_1f\|_{L^p_x} \lesssim \eta^{3(\frac1r-\frac1p)}\|f\|_{L^r_x},
\end{align}
and also by scaling argument, one has 
\begin{align}
\|B^\omega_\eta P_\lambda f\|_{L^p_x} \lesssim \eta^{3(\frac1r-\frac1p)}\lambda^{4(\frac1r-\frac1p)}\|f\|_{L^r_x}.
\end{align}
\end{lem}
We observe that in Lemma \ref{sobolev-emb}, the function $f$ is not required to be localised in an angular sector. Thus we can use Lemma \ref{sobolev-emb} first and then apply Lemma \ref{ang-con-est}. For example, we write
\begin{align*}
\left\|\sup_{\omega}\|S^\omega_{\mu,d}u\|_{L^\infty_x}\right\|_{L^2_t} & \lesssim 	\mu^{\frac43-}\left(\frac d\mu\right)^{1-}\left\|\sup_\omega\|S^\omega_{\mu,d}u\|_{L^{3+}_x}\right\|_{L^2_t} \\
& \lesssim \mu^{\frac43-}\left(\frac d\mu\right)^{\frac54-}\|S_{\mu,d}\langle\Omega\rangle^\frac12u\|_{L^2_tL^{3+}_x} \lesssim \mu^{\frac43-}\left(\frac d\mu\right)^{\frac54-}\|S_\mu u\|_{F_{\Omega,\mu}}.
\end{align*}

\section{Bilinear decompositions for angles}

We shall discuss various bilinear decompositions for frequency localised products of the form:
\begin{equation}\label{loc-prod}
S_{\lambda_0}(S_{\lambda_1}u\, S_{\lambda_2}v).
\end{equation}
By the standard Littlewood-Paley trichotomy, the localised product \eqref{loc-prod} vanishes unless
$$
\min\{\lambda_0,\lambda_1,\lambda_2 \} \lesssim \textrm{med}\{ \lambda_0,\lambda_1,\lambda_2 \} \approx \max\{\lambda_0,\lambda_1,\lambda_2\}.
$$
We focus on two important interactions in \eqref{loc-prod}, namely, the High$\times$High and Low$\times$High frequency interactions: 
\begin{align}
\lambda_0\lesssim & \, \lambda_1\approx\lambda_2, \label{hh} \\
\lambda_1\lesssim &  \, \lambda_0\approx\lambda_2. \label{lh}
\end{align}
In what follows, we give the bilinear decompositions for angles for the type \eqref{hh}, \eqref{lh}. We refer the reader to \cite{sterbenz,sterbenz2,wang} for more details. 
\begin{lem}[Lemma 6.1 of \cite{sterbenz2}]\label{hh-ang-decom}
For the following localised products of the form:
$$
S_{\mu,d}(S_{\lambda,\cdot\le\min\{d,c\mu\}}u\, S_{\lambda,\cdot\le\min\{d,c\mu\}}v),
$$
one has the following angular decomposition:
\begin{align*}
s_{\mu,d}^\pm(s^-_{\lambda,\cdot\le\min\{d,c\mu\}}*s^+_{\lambda,\cdot\le\min\{d,c\mu\}}) & = \sum_{\substack{ \omega_1,\omega_3 \\ |\omega_1\mp\omega_3|\approx(\frac{d}{\mu})^\frac12 }}s_{\mu,d}^{\omega_1,\pm}(s^-_{\lambda,\cdot\le\min\{d,c\mu\}}*b^{\omega_3}_{(\frac{d}{\mu})^\frac12}s^+_{\lambda,\cdot\le\min\{d,c\mu\}}) \\
& = \sum_{\substack{ \omega_1,\omega_2,\omega_3 \\ |\omega_1\mp\omega_3|\approx(\frac{d}{\mu})^\frac12 \\ |\omega_2+\omega_3|\approx (\frac{d}{\mu})^\frac12 }}s_{\mu,d}^{\omega_1,\pm}(b^{\omega_2}_{(\frac{d}{\mu})^\frac12}s^-_{\lambda,\cdot\le\min\{d,c\mu\}}*b^{\omega_3}_{(\frac{d}{\mu})^\frac12}s^+_{\lambda,\cdot\le\min\{d,c\mu\}}),
\end{align*}
for the convolution of the associated cutoff functions in Fourier side. There is a similar decomposition for the terms $S_{\mu,\cdot\le d}(S_{\lambda,\cdot\le d}u\, S_{\lambda,d}v)$ and $S_{\mu,\cdot\le d}(S_{\lambda,d}u\, S_{\lambda,\cdot\le d}v)$, where $d$ is in the range $d<c\mu$ and $c\ll1$ is the small number fixed above.
\end{lem}
We also have the angular decomposition for the interaction of type \eqref{lh}. 
\begin{lem}[Lemma 6.2 of \cite{sterbenz2}]\label{lh-ang-decom}
For the following localised products of the form:
$$
S_{\lambda,\cdot\le\min\{d,c\mu\}}(S_{\mu,d}u\, S_{\lambda,\cdot\le\min\{d,c\mu\}}v),
$$
one has the following angular decomposition:
\begin{align*}
s^+_{\lambda,\cdot\le\min\{d,c\mu\}}(s^\pm_{\mu,d}*s^+_{\lambda,\cdot\le\min\{d,c\mu\}})  = \sum_{\substack{ \omega_1,\omega_2,\omega_3 \\ |\omega_1\mp\omega_2|\approx(\frac d\mu)^\frac12 \\ |\omega_1-\omega_3|\approx (\frac d\mu)^\frac12 }} b^{\omega_1}_{(\frac d\mu)^\frac12}s^+_{\lambda,\cdot\le\min\{d,c\mu\}}(s^{\omega_2,\pm}_{\mu,d}*b^{\omega_3}_{(\frac d\mu)^\frac12}s^+_{\lambda,\cdot\le\min\{d,c\mu\}}),
\end{align*}
for the convolution of the associated cutoff functions in Fourier side. There is a similar decomposition for the terms $S_{\lambda,\cdot\le d}(S_{\mu,\cdot\le d}u\, S_{\lambda,d}v)$ and $S_{\lambda,d}(S_{\mu,\cdot\le d}u\, S_{\lambda,\cdot\le d}v)$ in the range $d<c\mu$, where $c\ll1$ is a fixed small number.
\end{lem}
We note that the angular sectors involved in the summation on the right-hand side of the decompositions in Lemma \ref{hh-ang-decom} and Lemma \ref{lh-ang-decom} are essentially diagonal. To avoid verbatim, we write angular decompositions throughout this paper as follows:
\begin{align}
\begin{aligned}
&S_{\mu,d}(S_{\lambda,\cdot\le\min\{d,c\mu\}}\varphi\, S_{\lambda,\cdot\le\min\{d,c\mu\}}\phi)\\
& \qquad\qquad= \sum_\omega S^{\pm\omega}_{\mu,d}(B^{-\omega}_{(\frac d\mu)^\frac12}S_{\lambda,\cdot\le\min\{d,c\mu\}}\varphi\, B^\omega_{(\frac d\mu)^\frac12}S_{\lambda,\cdot\le\min\{d,c\mu\}}\phi), \label{hh-ang-decomp} 
\end{aligned}
\end{align}
\begin{align}
\begin{aligned}
&S_{\lambda,\cdot\le\min\{d,c\mu\}}(S_{\mu,d}\varphi\, S_{\lambda,\cdot\le\min\{d,c\mu\}}\phi) \\
&\qquad\qquad = \sum_\omega B^{\omega}_{(\frac d\mu)^\frac12}S_{\lambda,\cdot\le\min\{d,c\mu\}}(S^{\pm\omega}_{\mu,d}\varphi\, B^{\omega}_{(\frac d\mu)^\frac12}S_{\lambda,\cdot\le\min\{d,c\mu\}}\phi). \label{lh-ang-decomp} 	
\end{aligned}
\end{align}
We also need the following angular decomposition, which will be only used in Section \ref{pf-lh-far-cone}.
\begin{lem}[Lemma 6.4 of \cite{sterbenz2}]\label{lh-ang-decomp-spe}
For the following expression:
$$
S^\omega_{\lambda,d}(S_\mu u\, S_{\lambda,\cdot<c\mu}v),
$$	
one has the following angular restriction:
$$
s^{\omega_1,+}_{\lambda,d}(s^\pm_\mu*s^+_{\lambda,\cdot<c\mu}) = s^{\omega_1,+}_{\lambda,d}(s^\pm_{\mu}*b^{\omega_3}_{(\frac{d}{c\lambda})^\frac12}s^+_{\lambda,\cdot<c\mu})
$$
for the convolution of the associated cutoff functions in Fourier space. Here the angles are restricted to the range $|\omega_1-\omega_3|\approx\left(\dfrac{d}{c\lambda}\right)^\frac12$.
\end{lem}

\section{Proof of main Theorem}
This section is devoted to the proof of our main result. The proof of scattering is followed by Theorem \ref{gwp}.  
See \cite{sterbenz,sterbenz2,wang} for details.
Now we focus on the proof of Theorem \ref{gwp}. In view of Duhamel's principle and Picard's iteration, we need to show the following nonlinear estimates: 
 \begin{align}
 \|VQ(\varphi,\phi)\|_{F_\Omega} & \lesssim \|\varphi\|_{F_\Omega}\|\phi\|_{F_\Omega}, \label{bi-phi} \\
 \|V(\phi_1\partial\phi_2)\|_{F_\Omega} & \lesssim \|\phi_1\|_{F_\Omega}\|\phi_2\|_{F_\Omega}, \label{bi-a} \\
 \|V(\phi\varphi\psi)\|_{F_\Omega}	& \lesssim \|\phi\|_{F_\Omega}\|\varphi\|_{F_\Omega}\|\psi\|_{F_\Omega}. \label{cubic}
 \end{align}
The bilinear estimates \eqref{bi-a} is already known. (See \cite{sterbenz2}.)
The treatment of cubic terms is very straightforward. It suffices to show that
\begin{align}
\|S_\mu(S_1\langle\Omega\rangle uS_1vS_1w)\|_{L^1_tL^2_x} & \lesssim \mu^{0+}\|S_1u\|_{F_\Omega}\|S_1v\|_{F_\Omega}\|S_1w\|_{F_\Omega}, \label{cubic-1} \\
\|S_1(S_\mu\langle\Omega\rangle uS_1vS_1w)\|_{L^1_tL^2_x} & \lesssim \mu\|S_\mu u\|_{F_{\Omega,\mu}}\|S_1v\|_{F_\Omega}\|S_1w\|_{F_\Omega}. \label{cubic-2}
\end{align}
We simply use the Bernstein's inequality, H\"{o}lder's inequality and the Strichartz estimates. Indeed, for the proof of \eqref{cubic-1}, we see that
\begin{align*}
\|S_\mu(S_1\langle\Omega\rangle uS_1vS_1w)\|_{L^1_tL^2_x} & \lesssim \mu^{2-}	\|S_\mu(S_1\langle\Omega\rangle uS_1vS_1w)\|_{L^1_tL^{1+}_x} \\
& \lesssim \mu^{2-}\|S_1\langle\Omega\rangle u\|_{L^\infty_tL^2_x}\|S_1v\|_{L^2_tL^6_x}\|S_1w\|_{L^{3+}_x} \\
& \lesssim \|S_1u\|_{F_\Omega}\|S_1v\|_{F_\Omega}\|S_1w\|_{F_\Omega}.
\end{align*}
To prove \eqref{cubic-2}, we write
\begin{align*}
	\|S_1(S_\mu\langle\Omega\rangle uS_1vS_1w)\|_{L^1_tL^2_x} & \lesssim \|S_\mu\langle\Omega\rangle u\|_{L^2_tL^\infty_x}\|S_1vS_1w\|_{L^2_tL^2_x} \\
	& \lesssim \mu^{\frac46}\|S_\mu\langle\Omega\rangle u\|_{L^2_tL^6_x}\|S_1v\|_{L^2_tL^\infty_x}\|S_1w\|_{L^\infty_tL^2_x} \\
	& \lesssim \mu^\frac32 \|S_\mu u\|_{F_{\Omega,\mu}}\|S_1v\|_{F_\Omega}\|S_1w\|_{F_\Omega},
\end{align*}
where we used $\|S_\mu\langle\Omega\rangle u\|_{L^2_tL^6_x}\lesssim \mu^\frac56\|u\|_{F_{\Omega,\mu}}$.

Fron now on we exclusively consider the bilinear estimates \eqref{bi-phi}. We first apply the dyadic decomposition on the space-time frequency of the bilinear form $Q(\varphi,\phi)$. It suffices to treat the High$\times$High and Low$\times$High interactions as follows:
\begin{align*}
\sum_{\mu\lesssim\max\{\lambda_1,\lambda_2\}}\mu \|VQ(S_{\lambda_1}\varphi,S_{\lambda_2}\phi)\|_{F_{\Omega,\mu}} & \lesssim \lambda_1\lambda_2\|\varphi\|_{F_{\Omega,\lambda_1}}\|\phi\|_{F_{\Omega,\lambda_2}}, \\
\|VQ(S_{\mu}\varphi,S_\lambda\phi)\|_{F_{\Omega,\lambda} } & \lesssim \mu\|\varphi\|_{F_{\Omega,\mu}}\|\phi\|_{F_{\Omega,\lambda}}.
\end{align*}
As we are concerned with a scale-invariant function space, it is reasonable to assume that the high frequency $\lambda=\lambda_1=\lambda_2=1$ and the low frequency $\mu\lesssim1$. In consequence, our aim is to prove the following bilinear estimates:
\begin{align}
\mu\|VQ(S_{1}\varphi,S_1\phi)\|_{F_{\Omega,\mu}} & \lesssim \mu^{0+}\|S_1\phi\|_{F_\Omega}\|S_1\phi\|_{F_\Omega},\label{bi-hh} \\
\|VQ(S_\mu\varphi,S_1\phi)\|_{F_{\Omega,1}} & \lesssim \mu\|S_\mu\varphi\|_{F_{\Omega,\mu}}\|S_1\phi\|_{F_\Omega}. \label{bi-lh}
\end{align}
We present the main scheme of the proof of \eqref{bi-hh} and \eqref{bi-lh}. We first note that the Fourier projection operator $S_\mu$ and the parametrix $V$ do not commute when $\mu\lesssim1$. Indeed, we have (see \cite{sterbenz2} for derivation,)
\begin{align}
	S_\mu VQ(S_1\varphi,S_1\phi) = S_\mu VS_\mu Q(S_1\varphi,S_1\phi) - \sum_{\mu\lesssim\sigma\lesssim1}W(P_\mu S_{\sigma,\sigma}VQ(S_1\varphi,S_1\phi)). \label{hh-commute}
\end{align}
 For the second term of the right-handside of \eqref{hh-commute}, we shall prove the following:
 \begin{align}
 \sum_\mu\mu\left\|\sum_{\mu\lesssim\sigma\lesssim1}P_\mu S_{\sigma,\sigma}VQ\langle\Omega\rangle(S_1\varphi,S_1\phi)\right\|_{L^\infty_tL^2_x} \lesssim \|S_1\varphi\|_{F_\Omega}\|S_1\phi\|_{F_\Omega},
 \end{align}
which is very straightforward. Then we further decompose $S_\mu Q(S_1\varphi,S_1\phi)$ into the space-time frequency which is away from the light cone and near the cone, respectively. For this purpose, we write
\begin{align}
\begin{aligned}
S_\mu Q(S_1\varphi,S_1\phi) & = S_\mu Q(S_1\varphi, S_{1,c\mu\le\cdot}\phi) + S_\mu Q(S_{1,c\mu\le\cdot}\varphi,S_{1,\cdot\le c\mu}\phi) \\
& \qquad\qquad	+S_\mu Q(S_{1,\cdot\le c\mu}\varphi,S_{1,\cdot\le c\mu}\phi)   \\
& \quad =: \mathcal H\mathcal H^1+\mathcal H\mathcal H^2+\mathcal H\mathcal H^3. 
\end{aligned}
\end{align}
The first and second terms $\mathcal H\mathcal H^1,\mathcal H\mathcal H^2$ are rather easier than the third term $\mathcal H\mathcal H^3$, which is near the cone. We will use the angular decomposition \eqref{hh-ang-decom} to $\mathcal H\mathcal H^3$. Then we apply the H\"{o}lder inequality, Sobolev embedding, and angular concentration estimates and the Strichartz estimates to obtain the required estimates. 

In the proof of \eqref{bi-lh}, we simply get
$$
S_1VQ(S_1\varphi,S_1\phi) = S_1VS_1Q(S_1\varphi,S_1\phi)
$$
and we write 
\begin{align}
\begin{aligned}
S_1Q(S_\mu\varphi, S_1\phi) & = S_1Q(S_\mu\varphi,S_{1,c\mu\le\cdot}\phi) + S_{1,c\mu\le\cdot}Q(S_\mu\varphi,S_{1,\cdot\le c\mu}\phi) \\
& \qquad\qquad + S_{1,\cdot\le c\mu}Q(S_\mu\varphi,S_{1,\cdot\le c\mu}\phi) 	 \\
& \quad =: \mathcal L\mathcal H^1+\mathcal L\mathcal H^2+\mathcal L\mathcal H^3. 
\end{aligned}
\end{align}
The Low$\times$High interaction is more difficult than the High$\times$High interaction. This is because of the regularity $\dot{B}^1$ of the scalar field $\varphi$. We will lose $\mu$ from $\varphi$, and hence we must gain more. To overcome this problem, we first observe that the angular decomposition as Lemma \ref{lh-ang-decom} will be very large when the modulation $d$ is close to the low frequency $\mu$. Thus we shall divide the Low$\times$High interaction into two cases: $(\frac d\mu)^\frac12\ll\mu$ and $\mu\lesssim(\frac d\mu)^\frac12$. We still use Lemma \ref{lh-ang-decom} for the first case. In the second case, instead of Lemma \ref{lh-ang-decom}, we use a smaller sector with size $\mu$. In this case, the range of $d$ is given by $\mu^3\lesssim d\lesssim \mu$ and hence the summation on $d$ makes no problem. We refer the readers to \cite{wang} for the change of weight between $\mu$ and $(\frac d\mu)^\frac12$. 

We introduce the outline of the remainder of this section. We first treat the proof of \eqref{hh-commute} in Section \ref{pf-hh-com}. Then the estimates of the High$\times$High and Low$\times$High interaction away from the light cone is given in Section \ref{pf-hh-far-cone} and \ref{pf-lh-far-cone}, which are rather easier than the frequency near the cone. Section \ref{pf-hh-near-cone} and \ref{pf-lh-near-cone} are devoted to the estimates of the High$\times$High and Low$\times$High interaction near the cone. 
\subsection{High$\times$High interaction including Commutator term}\label{pf-hh-com}
We simply use the Bernstein's inequality, the boundedness of multipliers (Lemma \ref{multi-bdd}), and H\"{o}lder's inequality.
\begin{align*}
\sum_{\mu\lesssim1}\mu \left\|\sum_{\mu\lesssim\sigma\lesssim1}P_\mu S_{\sigma,\sigma}VQ(S_1\varphi,S_1\phi)\right\|_{L^\infty_tL^2_x} & \lesssim \sum_{\mu\lesssim1}\mu\sum_{\mu\lesssim\sigma\lesssim1}\|P_\mu S_{\sigma,\sigma}VQ(S_1\varphi,S_1\phi)\|_{L^\infty_tL^2_x} \\
& \lesssim \sum_{\mu\lesssim1}\mu\sum_{\mu\lesssim\sigma\lesssim1}\mu^2\|S_{\sigma,\sigma}VQ(S_1\varphi,S_1\phi)\|_{L^\infty_tL^1_x} \\
& \lesssim \sum_{\mu\lesssim1}\mu\sum_{\mu\lesssim\sigma\lesssim1}\left(\frac\mu\sigma\right)^2\|S_1\varphi S_1\phi\|_{L^\infty_tL^1_x} \\
& \lesssim \sum_{\mu\lesssim1}\mu\|S_1\varphi\|_{L^\infty_tL^2_x}\|S_1\phi\|_{L^\infty_tL^2_x}	\\
& \lesssim \|S_1\varphi\|_{F_\Omega}\|S_1\phi\|_{F_\Omega}.
\end{align*}


\subsection{High$\times$High interaction away from cone }\label{pf-hh-far-cone}
We recall the definition of the $Y$ space. We first observe that
$$
\sum_\mu \|S_\mu (S_1\varphi S_1\phi)\|_{L^1_tL^2_x} \lesssim \sum_{\mu\lesssim1}\sum_{\sigma\lesssim\mu}\|P_\sigma S_\mu (S_1\varphi S_1\phi)\|_{L^1_tL^2_x}.
$$
We use the Bernstein's inequality to gain $\sigma^\frac23$. The H\"{o}lder inequality and Strichartz estimates give the desired estimates as follows. 
\begin{align*}
\mu\|V\mathcal H\mathcal H^1\|_{F_{\Omega,\mu}}  \lesssim \mu\|V\mathcal H\mathcal H^1\|_{\langle\Omega\rangle^{-1}Y_{\mu}} 
& \lesssim \sum_{\mu\lesssim1}\sum_{\sigma\lesssim\mu}\sigma^\frac23\|\langle\Omega\rangle(S_1\varphi S_{1,c\mu\le\cdot}\phi)\|_{L^1_tL^\frac32_x} \\
& \lesssim \sum_{\mu\lesssim1}\sum_{\sigma\lesssim\mu}\sigma^\frac23\|S_1\langle\Omega\rangle\varphi\|_{L^2_tL^6_x}\|S_{1,c\mu\le\cdot}\langle\Omega\rangle\phi\|_{L^2_tL^2_x} \\
& \lesssim \sum_{\mu\lesssim1}\sum_{\sigma\lesssim\mu}\sigma^\frac23 (c\mu)^{-\frac12}\|S_1\varphi\|_{F_\Omega}\|S_1\phi\|_{F_\Omega} \\
& \lesssim \sum_{\mu\lesssim1}c^{-\frac12}\mu^\frac16\|S_1\varphi\|_{F_\Omega}\|S_1\phi\|_{F_\Omega} \\
& \lesssim c^{-\frac12}\|S_1\varphi\|_{F_\Omega}\|S_1\phi\|_{F_\Omega}, 
\end{align*}
where we used $\mu^\frac12\|S_{1,\mu\le\cdot}\langle\Omega\rangle u\|_{L^2_tL^2_x} \lesssim \|S_1u\|_{F_{\Omega}}$.
The estimate of $\mathcal H\mathcal H^2$ is very similar. Indeed,
\begin{align*}
\mu\|V	\mathcal H\mathcal H^2\|_{F_{\Omega,\mu}} & \lesssim \sum_{\mu\lesssim1}\sum_{\sigma\lesssim\mu}\|P_\sigma S_\mu Q\langle\Omega\rangle(S_{1,c\mu\le\cdot}\varphi,S_{1,\cdot\le c\mu}\phi)\|_{L^1_tL^2_x} \\
& \lesssim \sum_{\mu\lesssim1}\sum_{\sigma\lesssim\mu}\sigma^\frac23\|\langle\Omega\rangle(S_{1,c\mu\le\cdot}\varphi S_{1,\cdot\le c\mu}\phi)\|_{L^1_tL^\frac32_x} \\
& \lesssim \sum_{\mu\lesssim1}\sum_{\sigma\lesssim\mu}\sigma^\frac23 (c\mu)^{-\frac12}\|S_1\varphi\|_{F_\Omega}\|S_1\phi\|_{F_\Omega} \\
& \lesssim c^{-\frac12}\mu^\frac16 \|S_1\varphi\|_{F_\Omega}\|S_1\phi\|_{F_\Omega}.
\end{align*}

\subsection{Low$\times$High interaction away from cone}\label{pf-lh-far-cone}
We first use the H\"{o}lder's inequality with respect to $t$ and then $x$ to get
\begin{align*}
\|V\mathcal L\mathcal H^1\|_{F_{\Omega,1}} \lesssim \|V\mathcal L\mathcal H^1\|_{\langle\Omega\rangle^{-1}Y_1} & = \|S_1Q\langle\Omega\rangle(S_\mu\varphi,S_{1,c\mu\le\cdot}\phi)\|_{L^1_tL^2_x}	 \\
& \lesssim \|S_\mu\langle\Omega\rangle\varphi\|_{L^2_tL^\infty_x}\|S_{1,c\mu\le\cdot}\langle\Omega\rangle\phi\|_{L^2_tL^2_x} \\
& \lesssim \mu^\frac23\|S_\mu\varphi\|_{L^2_tL^6_x}\|S_{1,c\mu\le\cdot}\phi\|_{L^2_tL^2_x} \\
& \lesssim \mu^\frac23\mu^\frac56 \|S_\mu\varphi\|_{F_{\Omega,\mu}}(c\mu)^{-\frac12}\|S_1\phi\|_{F_\Omega} \\
& \lesssim c^{-\frac12}\mu \|S_\mu\varphi\|_{F_{\Omega,\mu}}\|S_1\phi\|_{F_\Omega}.
\end{align*}
Here we also used $\|S_\mu\langle\Omega\rangle\varphi\|_{L^2_tL^6_x}\lesssim\mu^\frac56\|S_\mu\varphi\|_{F_{\Omega,\mu}}$. The estimate of $\mathcal L\mathcal H^2$ in the $X$ space is quite similar.
\begin{align*}
\|V\mathcal L\mathcal H^2\|_{F_\Omega} \lesssim \|V\mathcal L\mathcal H^2\|_{\langle\Omega\rangle^{-1}X^{\frac12,1}_1} & \lesssim \sum_{c\mu<d}d^{-\frac12}\|S_{1,d}Q\langle\Omega\rangle(S_\mu\varphi, S_{1,\cdot\le c\mu}\phi)\|_{L^2_tL^2_x} \\
& \lesssim (c\mu)^{-\frac12}\|S_\mu\langle\Omega\rangle\varphi\|_{L^2_tL^\infty_x}\|S_{1,\cdot\le c\mu}\langle\Omega\rangle\phi\|_{L^\infty_tL^2_x} \\
& \lesssim c^{-\frac12}\mu \|S_\mu\varphi\|_{F_{\Omega,\mu}}\|S_1\phi\|_{F_\Omega}.	
\end{align*}
We need to estimate $\mathcal L\mathcal H^2$ in the $Z$ space. To do this, we apply Lemma \ref{lh-ang-decomp-spe}, H\"{o}lder inequality with respect to $t$ and the Bernstein's inequality and then Lemma \ref{ang-con-est}.
\begin{align*}
\|V\mathcal L\mathcal H^2\|_{Z_{\Omega,1}} & = \sum_{c\mu<d}d^{-1}\int\sup_\omega\|B^\omega_{d^\frac12}S_{1,d}Q(S_\mu, S_{1,\cdot\le c\mu}\phi)\|_{L^\infty_x}dt \\
& \lesssim \sum_{c\mu<d}d^{-1}\int\sup_{\substack{\omega_1,\omega_2 \\ |\omega_1-\omega_2|\approx(\frac dc)^\frac12}}\|B^{\omega_1}_{d^\frac12}S_{1,d}Q(S_\mu\varphi,B^{\omega_2}_{(\frac dc)^\frac12}S_{1,\cdot\le c\mu}\phi)\|_{L^\infty_x}dt	\\
& \lesssim \sum_{c\mu<d}d^{-1}\|S_\mu\varphi\|_{L^2_tL^\infty_x}\left\|\sup_\omega\|B^\omega_{(\frac dc)^\frac12}S_{1,\cdot\le c\mu}\phi\|_{L^\infty_x}\right\|_{L^2_t} \\
& \lesssim \sum_{c\mu<d}d^{-\frac14-}\mu^{\frac43-}\|S_\mu\varphi\|_{L^2_tL^{3+}_x}\|S_1\langle\Omega\rangle^{\frac12}\phi\|_{L^2_tL^{3+}_x} \\
& \lesssim \mu^{\frac{13}{12}-}\|S_\mu\varphi\|_{F_{\Omega,\mu}}\|S_1\phi\|_{F_\Omega}.
\end{align*}

\subsection{High$\times$High interaction near cone}\label{pf-hh-near-cone}
We further decompose $\mathcal H\mathcal H^3$ as follows:
\begin{align*}
S_\mu Q(S_{1,\cdot\le c\mu}\varphi,S_{1,\cdot\le c\mu}\phi) & = \sum_{d<c\mu}S_{\mu,\cdot\le d}Q(S_{1,\cdot\le d}\varphi,S_{1,d}\phi)+ \sum_{d<c\mu}S_{\mu,\cdot\le d}Q(S_{1,d}\varphi,S_{1,\cdot\le d}\phi) \\
&\quad + \sum_{d\lesssim\mu}S_{\mu,d}Q(S_{1,\cdot\le\min\{d,c\mu\}}\varphi,S_{1,\cdot\le\min\{d,c\mu\}}\phi) \\
& =: \mathcal H\mathcal H^3_1+\mathcal H\mathcal H^3_2+\mathcal H\mathcal H^3_3.	
\end{align*}
We estimate the term $\mathcal H\mathcal H^3_1$ in the $Y$ space. We apply in order the angular decomposition (Lemma \ref{hh-ang-decom}), Sobolev embedding estimate (Lemma \ref{sobolev-emb}), H\"{o}lder inequality, and then Strichartz estimates.
\begin{align*}
\mu\|V\mathcal H\mathcal H^3_1\|_{\langle\Omega\rangle^{-1}Y_{\mu}} 
& \lesssim \sum_{d<c\mu}\left\| \left(\sum_{\omega}\|S^{\omega}_{\mu,\cdot\le d}\langle\Omega\rangle Q(B^{-\omega}_{(\frac d\mu)^\frac12}S_{1,\cdot\le d}\varphi,B^\omega_{(\frac d\mu)^\frac12}S_{1,d}\varphi)\|_{L^2_x}^2	\right)^\frac12\right\|_{L^1_t} \\
& \lesssim \sum_{d<c\mu}\mu^\frac23 \left(\frac d\mu\right)^\frac14\left(\frac d\mu\right)^\frac12 \left\| \left(\sum_{\omega}\|\langle\Omega\rangle (B^{-\omega}_{(\frac d\mu)^\frac12}S_{1,\cdot\le d}\varphi,B^\omega_{(\frac d\mu)^\frac12}S_{1,d}\varphi)\|_{L^\frac32_x}^2	\right)^\frac12\right\|_{L^1_t} \\
& \lesssim \sum_{d<c\mu}\mu^{-\frac{1}{12}}d^\frac34\left\| \sup_\omega \|B^{-\omega}_{(\frac d\mu)^\frac12}S_{1,\cdot\le d}\langle\Omega\rangle\varphi\|_{L^6_x} \left(\sum_\omega\|B^\omega_{(\frac d\mu)^\frac12}S_{1,d}\langle\Omega\rangle\phi\|_{L^2_x}^2\right)^\frac12\right\|_{L^1_t} \\
& \lesssim \sum_{d<c\mu}\mu^{-\frac{1}{12}}d^\frac34\|S_{1,\cdot\le d}\langle\Omega\rangle\varphi\|_{L^2_tL^6_x}\|S_{1,d}\langle\Omega\rangle\phi\|_{L^2_tL^2_x} \\
& \lesssim \mu^\frac16 \|S_1\varphi\|_{F_\Omega}\|S_1\phi\|_{F_\Omega}.
\end{align*}
The estimate of $\mathcal H\mathcal H^3_2$ is very similar. We omit the details. The treatment of $\mathcal H\mathcal H^3_3$ is also similar. Indeed, we apply the angular decomposition, Sobolev estimates and H\"{o}lder inequality.
\begin{align*}
\mu\|V\mathcal H\mathcal H^3_3\|_{\langle\Omega\rangle^{-1}X^{\frac12,1}_\mu} & = \sum_{d\lesssim\mu}d^{-\frac12}\|S_{\mu,d}\langle\Omega\rangle Q(S_{1,\cdot\le\min\{d,c\mu\}}\varphi,S_{1,\cdot\le\min\{d,c\mu\}}\phi)\|_{L^2_tL^2_x} \\
& \lesssim \sum_{d\lesssim\mu}d^{-\frac12}\mu^{\frac{5}{12}}d^\frac14  \left(\frac d\mu\right)^\frac12     \\ 
&\qquad\times\left\| \left(\sum_\omega\|\langle\Omega\rangle(B^{-\omega}_{(\frac d\mu)^\frac12}S_{1,\cdot\le\min\{d,c\mu\}}\varphi S_{1,\cdot\le\min\{d,c\mu\}}\phi)\|^2_{L^\frac32_x}\right)^\frac12\right\|_{L^2_t} \\
& \lesssim \sum_{d\lesssim\mu}\mu^{-\frac{1}{12}}d^\frac14\|S_1\langle\Omega\rangle\varphi\|_{L^2_tL^6_x}\|S_1\langle\Omega\rangle\phi\|_{L^\infty_tL^2_x} \\
& \lesssim \mu^\frac16 \|S_1\varphi\|_{F_\Omega}\|S_1\phi\|_{F_\Omega}.
\end{align*}
We need to estimate the term $\mathcal H\mathcal H^3_3$ in the $Z$ space also. We first use Lemma \ref{multi-bdd} and then the following step is quite similar.
\begin{align*}
\mu\|V\mathcal H\mathcal H^3_3\|_{Z_{\Omega,\mu}} 
& \lesssim \sum_{d\lesssim\mu}(d\mu)^{-1}\int\sup_\omega\|S^\omega_{\mu,d}Q(S_{1,\cdot\le\min\{d,c\mu\}}\varphi,S_{1,\cdot\le\min\{d,c\mu\}}\phi)\|_{L^\infty_x}dt	\\
& \lesssim \sum_{d\lesssim\mu}(d\mu)^{-1}\mu^{\frac83-}\left(\frac d\mu\right)^{1-}\left(\frac d\mu\right)^{\frac12}\int \| S_{1,\cdot\le\min\{d,c\mu\}}\varphi,S_{1,\cdot\le\min\{d,c\mu\}}\phi\|_{L^{\frac32+}_x}dt \\
& \lesssim \sum_{d\lesssim\mu}\mu^{\frac16-}d^{\frac12-}\int \|S_{1,\cdot\le\min\{d,c\mu\}}\varphi\|_{L^{3+}_x}\|S_{1,\cdot\le\min\{d,c\mu\}}\phi\|_{L^{3+}_x}dt \\
& \lesssim \mu^{\frac23-}\|S_1\varphi\|_{L^2_tL^{3+}_x}\|S_1\phi\|_{L^2_tL^{3+}_x} \\
& \lesssim \mu^{\frac23-}\|S_1\varphi\|_{F_\Omega}\|S_1\phi\|_{F_\Omega}.
\end{align*}

\subsection{Low$\times$High interaction near cone}\label{pf-lh-near-cone}


As the previous section, we further decompose $\mathcal L\mathcal H^3$ as follows:
\begin{align*}
S_{1,\cdot\le c\mu}Q(S_{\mu}\varphi,S_{1,\cdot<c\mu}\phi) & = \sum_{d<c\mu}S_{1,\cdot\le d}Q(S_{\mu,\cdot\le d}\varphi,S_{1,d}\phi)+ \sum_{d<c\mu}S_{1,d}Q(S_{\mu,\cdot\le d}\varphi,S_{1,\cdot\le d}\phi) \\
& \quad + \sum_{d\lesssim\mu}S_{1,\cdot\le\min\{d,c\mu\}}Q(S_{\mu,d}\varphi,S_{1,\cdot\le\min\{d,c\mu\}}\phi)\\
& =: \mathcal L\mathcal H^3_1+\mathcal L\mathcal H^3_2+\mathcal L\mathcal H^3_3.	
\end{align*}
To deal with the term $\mathcal L\mathcal H^3_1$ and $\mathcal L\mathcal H^3_2$, we first apply Lemma \ref{lh-ang-decom}. Then we use in order the H\"{o}lder inequality with respect to $t$, and then $x$. Then we use Lemma \ref{sobolev-emb} and Strichartz estimates to obtain the required estimates. 
The explicit treatment of $\mathcal L\mathcal H^3_1,\mathcal L\mathcal H^3_2$ is as follows.
\begin{align*}
\|V\mathcal L\mathcal H^3_1\|_{\langle\Omega\rangle^{-1}Y_1} 
& \lesssim \sum_{d<c\mu}\left(\frac d\mu\right)^\frac12 \left\|\sup_\omega\|S^\omega_{\mu,\cdot\le d}\langle\Omega\rangle\varphi\|_{L^\infty_x}\right\|_{L^2_t}\left\| \left(\sum_\omega\|B^\omega_{(\frac d\mu)^\frac12}S_{1,d}\langle\Omega\rangle\phi\|_{L^2_x}^2\right)^\frac12\right\|_{L^2_t} \\
& \lesssim \sum_{d<c\mu}\mu^{-\frac{1}{12}}d^\frac34 \|S_{\mu}\langle\Omega\rangle\varphi\|_{L^2_tL^6_x}d^{-\frac12}\|S_1\phi\|_{F_\Omega} \\
& \lesssim \sum_{d<c\mu}\mu^{-\frac{1}{12}}d^\frac14\mu^\frac56\|S_\mu\varphi\|_{F_{\Omega,\mu}}\|S_1\phi\|_{F_\Omega} \\
& \lesssim \mu \|S_\mu\varphi\|_{F_{\Omega,\mu}}\|S_1\phi\|_{F_\Omega},
\end{align*}
\begin{align*}
\|V\mathcal L\mathcal H^3_2\|_{F_{\Omega,1}}  \lesssim \|V\mathcal L\mathcal H^3_2\|_{\langle\Omega\rangle^{-1}X^{\frac12,1}_1}
& \lesssim \sum_{d<c\mu}d^{-\frac12}\left(\frac d\mu\right)^\frac12 \\
&\qquad\times \left\|\sup_\omega \|S^\omega_{\mu,\cdot\le d}\langle\Omega\rangle\varphi\|_{L^\infty_x} \left(\sum_\omega\|B^\omega_{(\frac d\mu)^\frac12}S_{1,\cdot\le d}\langle\Omega\rangle\phi\|_{L^2_x}^2\right)^\frac12\right\|_{L^2_t} \\
& \lesssim \sum_{d<c\mu}\mu^\frac16\left(\frac d\mu\right)^\frac14 \|S_{\mu}\langle\Omega\rangle\varphi\|_{L^2_tL^6_x}\|S_1\langle\Omega\rangle\phi\|_{L^\infty_tL^2_x} \\
& \lesssim \mu\|S_\mu\varphi\|_{F_{\Omega,\mu}}\|S_1\phi\|_{F_\Omega}.  
\end{align*}
We also need to estimate the term $\mathcal L\mathcal H^3_2$ in the $Z$ space. Fron now on, the proof is quite different. 
\begin{align*}
\|V\mathcal L\mathcal H^3_2\|_{Z_{\Omega,1}} & = \sum_{d\lesssim\mu}\int\sup_\omega\|B^\omega_{(\frac d\mu)^\frac12}S_{1,d}V\mathcal L\mathcal H^3_2\|_{L^\infty_x}dt	\\
& \lesssim \sum_{d\lesssim\mu}d^{-1}\int\sup_\omega\|B^\omega_{(\frac d\mu)^\frac12}S_{1,d}Q(S_{\mu,\cdot\le d}\varphi,S_{1,\cdot\le d}\phi)\|_{L^\infty_x}dt \\
& \lesssim \sum_{d\lesssim\mu}(d\mu)^{-\frac12}\int \sup_\omega\|B^\omega_{(\frac d\mu)^\frac12}S_{1,d}(S^\omega_{\mu,\cdot\le d}\varphi B^\omega_{(\frac d\mu)^\frac12}S_{1,\cdot\le d})\|_{L^\infty_x}dt \\
& \lesssim \sum_{d\lesssim\mu}(d\mu)^{-\frac12}\int\sup_\omega\|S^\omega_{\mu,\cdot\le d}\varphi\|_{L^\infty_x}\sup_\omega\|B^\omega_{(\frac d\mu)^\frac12}S_{1,\cdot\le d}\phi\|_{L^\infty_x}dt \\
& \lesssim \sum_{d\lesssim\mu}(d\mu)^{-\frac12} \mu^{\frac43-}\left(\frac d\mu\right)^{1-} \int\|S_{\mu}\varphi\|_{L^{3+}_x}\|S_1\phi\|_{L^{3+}_x}dt \\
&\lesssim \sum_{d\lesssim\mu}\mu^{-\frac16}d^{\frac12-}\|S_\mu\varphi\|_{L^2_tL^{3+}_x}\|S_1\phi\|_{L^2_tL^{3+}_x}. 
\end{align*}
In this manner, we can obtain the required estimate only if $d\ll\mu^3$. On the other hand, for $\mu^3\lesssim d<c\mu$, we have $\mu\lesssim(\frac d\mu)^\frac12$, and hence instead of Lemma \ref{lh-ang-decom}, we make the use of a smaller angular decomposition with size $\mu$ for output frequency and high input frequency. 
\begin{align*}
& \sum_{\mu^3\lesssim d<c\mu}(d\mu)^{-\frac12}\int\sup_\omega\|B^{\omega}_{(\frac d\mu)^\frac12}S_{1,d}(S_{\mu,\cdot\le d}\varphi S_{1,\cdot\le d}\phi)\|_{L^\infty_x}dt   \\
& \lesssim  \sum_{\mu^3\lesssim d<c\mu}(d\mu)^{-\frac12}\int\sup_\omega\|B^\omega_\mu S_{1,d}(S_{\mu,\cdot\le d}\varphi S_{1,\cdot\le d}\phi)\|_{L^\infty_x}dt \\
 & \lesssim \sum_{\mu^3\lesssim d<c\mu}(d\mu)^{-\frac12}\int \sup_{\omega_1}\sum_{\substack{|\omega_1+\omega_3|\approx \mu \\ |\omega_2\mp\omega_3|\approx(\frac d\mu)^\frac12}}\|B^{\omega_1}_\mu S_{1,d}(S^{\omega_2}_{\mu,\cdot\le d} B^{\omega_3}_\mu S_{1,\cdot\le d}\phi )\|_{L^\infty_x}dt \\
& \lesssim 	\sum_{\mu^3\lesssim d<c\mu }\sum_{\omega_2,\omega_3}(d\mu)^{-\frac12}\int \|S^{\omega_2}_{\mu,\cdot\le d}\varphi\|_{L^\infty_x}\|B^{\omega_3}_\mu S_{1,\cdot\le d}\phi\|_{L^\infty_x}dt \\
& \lesssim \sum_{\mu^3\lesssim d<c\mu}\sum_{\omega_2,\omega_3}(d\mu)^{-\frac12}\mu^{\frac43-}\left(\frac d\mu\right)^{\frac12-}\mu^{1-} \|S^{\omega_2}_{\mu,\cdot\le d}\varphi\|_{L^2_tL^{3+}_x}\|B^{\omega_3}_\mu S_{1,\cdot\le d}\phi\|_{L^2_tL^{3+}_x} \\
& \lesssim \mu^{\frac43-}\|S_\mu\varphi\|_{F_{\Omega,\mu}}\|S_1\phi\|_{F_\Omega}.
\end{align*}
Here, the summation on $d$ makes only the loss of $\log\mu$ and hence we get the desired estiates. 
The term $\mathcal L\mathcal H^3_3$ is the most crucial part of our proof. We write
\begin{align*}
\|V\mathcal L\mathcal H^3_3\|_{F_{\Omega,1}} & \lesssim \sum_{d\lesssim\mu}\|S_{1,\cdot\le\min\{d,c\mu\}}Q(S_{\mu,d}\varphi,S_{1,\cdot\le\min\{d,c\mu\}}\langle\Omega\rangle\phi)\|_{L^1_tL^2_x} \\
& \quad + \sum_{d\lesssim\mu}\|S_{1,\cdot\le\min\{d,c\mu\}}Q(S_{\mu,d}\langle\Omega\rangle\varphi,S_{1,\cdot\le\min\{d,c\mu\}}\phi)\|_{L^1_tL^2_x} \\
& =: \mathcal J_{1}+\mathcal J_{2}.	
\end{align*}
The term $\mathcal J_1$ is rather easier than $\mathcal J_2$. We simply recall the property of the $Z$ space \eqref{prop-z}.
\begin{align*}
\mathcal J_{1} & \lesssim \sum_{d\lesssim\mu}\left(\frac d\mu\right)^\frac12\left\|\sup_\omega\|S^\omega_{\mu,d}\varphi\|_{L^\infty_x}\right\|_{L^1_t}\|S_{1,\cdot\le\min\{d,c\mu\}}\langle\Omega\rangle\phi\|_{L^\infty_tL^2_x} \\
& \lesssim \mu\|S_\mu\varphi\|_{F_{\Omega,\mu}}\|S_1\phi\|_{F_\Omega}.
\end{align*}
Now we are left to consider the $\mathcal J_2$. We further decompose the range of $d $ into $d\ll\mu^3$ and $\mu^3\lesssim d\lesssim\mu$. If $d\ll\mu^3$, then we apply Lemma \ref{lh-ang-decom}, H\"{o}lder inequality and then Sobolev estimates.
\begin{align*}
\mathcal J_{2}^{d\ll\mu^3} & \lesssim \sum_{d\ll\mu^3}\left(\frac d\mu\right)^\frac12 \left\|\left(\sum_\omega\|S^\omega_{\mu,d}\langle\Omega\rangle\varphi\|_{L^4_x}^2\right)^\frac12\right\|_{L^2_t}\left\|\sup_\omega\|B^\omega_{(\frac d\mu)^\frac12}S_{1,\cdot\le\min\{d,c\mu\}}\phi\|_{L^4_x}\right\|_{L^2_t} \\
& \lesssim \sum_{d\ll\mu^3}\left(\frac d\mu\right)^\frac12\mu\left(\frac d\mu\right)^\frac38\|S_{\mu,d}\langle\Omega\rangle\varphi\|_{L^2_tL^2_x}\left(\frac d\mu\right)^{\frac18-}\left\|\sup_\omega\|B^\omega_{(\frac d\mu)^\frac12}S_{1,\cdot\le\min\{d,c\mu\}}\phi\|_{L^{3+}_x}\right\|_{L^2_t} \\
& \lesssim \sum_{d\ll\mu^3}\mu^\frac12\left(\frac d\mu\right)^{\frac12-}\|S_\mu\varphi\|_{F_{\Omega,\mu}}\|S_1\phi\|_{L^2_tL^{3+}_x}.
\end{align*}
Thus we can get the desired estimate for $d\ll\mu^3$. On the other hand, if $\mu^3\lesssim d\lesssim\mu$, we write
\begin{align*}
\mathcal J_{2}^{\mu^3\lesssim d\lesssim\mu} = \sum_{d\lesssim\mu}\left\|\sum_{\substack{|\omega_1+\omega_3|\approx\mu \\ |\omega_2\mp\omega_3|\approx(\frac d\mu)^\frac12}}B^{\omega_1}_\mu S_{1,\cdot\le\min\{d,c\mu\}}Q(S^{\omega_2}_{\mu,d}\langle\Omega\rangle\varphi,B^{\omega_3}_\mu S_{1,\cdot\le\min\{d,c\mu\}}\phi)\right\|_{L^1_tL^2_x}.
\end{align*}
We rewrite the $L^2_x$ norm via duality and then use H\"{o}lder inequality.
\begin{align*}
\|\cdots\|_{L^2_x} \lesssim \sup_{\|h\|_{L^2_x}=1}\left(\frac d\mu\right)^\frac12\sum_{\omega_2}\|S^{\omega_2}_{\mu,d}\langle\Omega\rangle\varphi\|_{L^2_x}\left\|\sum_{\omega_1,\omega_2,\omega_3}B^{\omega_2}_{(\frac d\mu)^\frac12}P_{\cdot\lesssim\mu}(B^{\omega_3}_\mu S_{1,\cdot\le\min\{d,c\mu\}}\phi B^{\omega_1}_\mu P_{\cdot\lesssim1}h)\right\|_{L^2_x}.
\end{align*}
Then
\begin{align*}
\mathcal J_{2}^{\mu^3\lesssim d\lesssim\mu} & \lesssim \sum_{d\lesssim\mu}\left(\frac d\mu\right)^\frac12d^{-\frac12}d^{\frac12}\|S_{\mu,d}\langle\Omega\rangle\varphi\|_{L^2_tL^2_x}\mu^{\frac43-}\left(\frac d\mu\right)^{\frac12-}\left\|\sup_\omega\|B^{\omega}_\mu S_{1,\cdot\le\min\{d,c\mu\}}\phi\|_{L^{3+}_x}\right\|_{L^2_t} \\
& \lesssim \sum_{d\lesssim\mu}\mu^{\frac56-}\left(\frac d\mu\right)^{\frac12-}\mu^\frac12\|S_1\langle\Omega\rangle^\frac12\phi\|_{L^2_tL^{3+}_x}\|S_\mu\varphi\|_{F_{\Omega,\mu}} \\
& \lesssim \mu^{\frac43-}\|S_\mu\varphi\|_{F_{\Omega,\mu}}\|S_1\phi\|_{F_\Omega}.	
\end{align*}
\begin{rem}\label{rem-lh-ang}
In the Low$\times$High regime, the most difficult interaction is when the low frequency controls the angular regularity. In this case, we cannot exploit the angular concentration estimates to gain some positive power of $\mu$. Hence we only use the Sobolev embedding estimates.	
\end{rem}

\section*{Acknowledgements}
The author is supported by NRF-2018R1D1A3B07047782 and NRF-2016K2A9A2A13003815.



\end{document}